\documentclass[12pt,reqno]{amsart}
\usepackage{hyperref}
\usepackage{soul}
\usepackage{ amssymb }
\numberwithin{equation}{section}
\usepackage{hyperref}
\usepackage{float}
\usepackage{pgf,tikz}

\newtheorem{theorem}{Theorem}[section]
\newtheorem{definition}[theorem]{Definition}

\newtheorem{proposition}[theorem]{Proposition}

\newtheorem{question}[theorem]{Question}

\newtheorem{conj}[theorem]{Conjecture}

\newtheorem{corollary}[theorem]{Corollary}

\newtheorem*{notation*}{Notation}

\begin{document}

\title{Recent results on homological properties of binomial edge ideal of graphs}

\author[P. Das]{Priya Das}

\address{Vellore Institute of Technology, Vellore, Tamilnadu-632014, India}
\email{priya.das@vit.ac.in}

\thanks{AMS Classification 2010. Primary: 13D02, 05E40, 16E05, 13D05, 16E65, 13F65}
\maketitle
\begin{abstract}
 In this article, we give a comprehensive survey of the recent progress of research on the binomial edge ideal of a graph since 2018. 
\end{abstract}

\section{Introduction} Let $G$ be a finite simple graph on $n$ vertices with set of edges, $E(G)$ with no isolated vertices. Let $S=K[x_1,x_2,\ldots,x_n,y_1,y_2,\ldots,y_n]$ be the polynomial ring over a field $\mathbb{K}$ and $f_{ij}=x_i y_j-x_j y_i, \{i,j\}\in E(G)$, the ideal $J_G$ in $S$ generated by $f_{ij}, 1\leq i<j \leq n$ is called the binomial edge ideal (BEI) of $G$. It was introduced by Herzog-Hibi-Kahle-Rauh in the contest of study of conditional independence ideals \cite{herzog_binomial} and independently by \cite{ohtani_binomial}. It can be also seen as a special case of the generalized binomial edge ideal as defined by \cite{GBEI_Rauh}, the collection of 2-minors of the generic $(2\times n)$-matrix.

In \cite{Madani_survey}, Madani has given an overview of Gr\"{o}bner bases, primary decomposition, and minimal graded free resolution of binomial edge ideal. The survey by Madani loc. cit. is a comprehensive survey of all the works till the year of its publication (2018), related to the binomial edge ideal in the context of homological studies. In this article, we will provide a survey of all the works mainly those that are related to the regularity of the binomial edge ideal, depth, Betti number, generalised binomial edge ideal, and parity binomial edge ideal, since \cite{Madani_survey}.

Gr\"{o}bner bases of binomial edge ideal are studied in \cite{herzog_binomial}, \cite{ohtani_binomial}, the authors of loc. cit. proved that the generators of the binomial edge ideal form a Gr\"{o}bner basis for closed graphs. In \cite{Koszul_ene}, the authors proved that $S/J_{G}$ is Koszul when $G$ is a closed graph. The primary decomposition of the binomial edge ideal is studied in \cite{herzog_binomial}, using the primary decomposition the authors have described the minimal prime ideals of $J_G$.

The study of the minimal free resolution of binomial edge ideals is one of the important topics in combinatorial commutative algebra. as the binomial ideals are geometrically and algebraically an important class. Also since the work of \cite{BPS} has a very thorough discussion for the monomial ideals, it makes sense to look for a combinatorial understanding of the situation in the binomial ideals. Among the binomial ideals, binomial edge ideals of a graph is a major stream of work currently. Also a number of works have been done in the binomial toric ideals associated to distributive lattices, namely the Hibi ideal \cite{ene}. The authors in \cite{PH} have given the generators of the first syzygy of the Hibi ideal of a planar lattice in terms of four types of sublattices.  In \cite{binomial_reg2}, the authors proved that BEI of a graph has a linear resolution if $G$ is a complete graph and vice versa.
 The question of a pure resolution of a BEI is discussed in \cite{Pure_resolution}.

One of the most important invariants of an ideal is the Castelnuovo-Mumford regularity. It is proved in \cite{reg_length}, for any graph $G$, $l\leq reg(S/J_{G})\leq n-1$, where $n$ is the number of vertices of $G$ and $l$ is the length of the largest induced path in $G$. Also, in \cite{reg_path}, authors proved that reg$(S/J_{G})=n-1$ if and only if $G$ is a path. Madani-Kiani in \cite{conjecture_pair} has given the following conjecture 
\begin{conj}\label{conj1}
Let $G$ be a simple graph and $c(G)$ the number of maximal cliques of $G$. Then
\begin{center}
 reg $(S/J_G)\leq c(G)$.
\end{center}
\end{conj}
Madani and Kiani proved the above conjecture for closed graphs and generalized block graphs in \cite{binomial_reg2}, and \cite{conjecture_generalized}. The conjecture in general remains as one of the important open problem in this area. For a block graph, $c(G)$ coincides with the number of blocks in $G$. Trees being a  subclass of the class of block graphs and combining it with above result we have that, for a tree $T$, which has $n$ vertices $c(T)=n-1$. In \cite{reg_block}, authors have shown that reg$(S/J_{G})\leq n-1$, thus proved the conjecture for some classes of trees.
Also, this conjecture is proved for the class of chordal graphs and for fan graphs of complete graphs in \cite{conjecture_fan}. Independently, the authors in loc. cit. prove this conjecture for chordal graphs in \cite{conjecture_chordal}. Recently, this conjecture has been proved for $P_4$-free graphs in \cite{conjecture_cograph}. For the classes of block graphs and semi block graph the conjecture has been settled in \cite{Arvind_Jahangir}. Recently, in \cite{Proof_conjecture}, authors proved the Conjecture \ref{conj1} for all graphs and has given a sharp bound of it. The authors proved that reg$(S/J_{G})\leq \eta{(G)}$ (for definition of $\eta{(G)}$, see \ref{defn}).

In the article \cite{CM_binomial}, the authors conjectured that the regularity of the binomial edge ideal of a graph coincides with the regularity of its initial ideal. If the conjecture holds for two graphs $G_1,G_2$ then it also holds for the join of two graphs, in Theorem 2.1 in \cite{binomial_reg3}.
The above result was used by the authors to completely characterize the BEI of regularity 3, (Theorem 3.2 \cite{binomial_reg3}). Threshold graphs which are also known as $(2K_2, C_4, P_4)$-free graphs are the graphs which have regularity 3. For more details about threshold graph, see \cite{threshold}. Whereas authors have characterized the binomial edge ideal of regularity 2 in \cite{binomial_reg2}. A new general upper bound for the regularity of $S/J_G$ has been given in terms of the dimension of a simplicial complex in \cite{reg_licci}. In \cite{reg_Gblock}, authors of loc. cit. obtained an improved lower and upper bound for the regularity of the BEI of a generalized block graph. Generalized block graphs are the generalization of block graphs and it is introduced in \cite{HS_join}. In \cite{construction_HibiBEI}, Hibi-Matsuda conjectured the following
\begin{conj}[\cite{construction_HibiBEI}, Conjecture 0.1] Let $G$ be a graph on $n$ vertices. Then \begin{center}
    reg$(S/J_G)\leq$ deg $h_{S/J_G}(t)$.
\end{center}
\end{conj}

 They proved that this conjecture is true if $S/J_{G}$ has a unique extremal Betti number. In \cite{Betti_bipartite}, Schenzel and Zafar have shown complete bipartite graphs have unique extremal Betti number. Also, Zafar and Zahid proved that $n$-cycle, $C_n$ has a unique extremal Betti number, see \cite{Betti_cycle}. In \cite{lower_block}, Herzog-Rinaldo characterized block graphs which admits unique extremal Betti number. In \cite{reg_Gblock}, the author characterized a generalized block graph which admits a unique extremal Betti number. In \cite{Arvind_Jahangir}, authors provide a sufficient condition for this conjecture.

The depth of an ideal is another commutative algebraic invariant of primary importance, in the context of BEI there have been a number of works in this area, but it is far from being complete. In \cite{graph_connectivity}, authors proved that depth of $S/J_G$ is $n+1$, for a connected block graph $G$. In \cite{Betti_cycle}, the authors showed that depth $S/J_G$ is $n,$ for $n>3$. In \cite{depth_GBG}, authors gave the depth of the BEI of generalised block graphs. For a connected graph $G$ on $n$-vertices the depth of binomial edge ideal is at most $n$ and further, if $G$ is a connected graph on $n$ vertices such that $S/J_G$ is Cohen-Macaulay then depth of it is exactly $n+1$, see, \cite{graph_connectivity}. Also, in the same article, the authors provided upper bounds for the depth of binomial edge ideal of a general graph in terms of certain graph invariants. Depth of binomial edge ideal for a cone graph remains invariant under the process of taking cone on connected graph, see \cite{depth_arvind}. In the same article they have given a formula for the depth of binomial edge ideal of join of two graphs. Recently, in \cite{malayeri_depth} they have completely characterized the graphs for which is depth of the BEI is 5.

Betti numbers of binomial edge ideal of cycle, complete bipartite graphs, trees, unicyclic graphs, and cone graph  are known, see \cite{Betti_cycle}, \cite{Betti_bipartite}, \cite{Betti_tree}, \cite{depth_arvind}. In \cite{Betti_tree}, the authors have given the generators the first syzygy of binomial edge ideal of cycle and unicyclic graphs explicitly. 

Generalized binomial edge ideal is introduced by Rauh in \cite{GBEI_Rauh} as a generalization to encompass a variety of binomial ideals. As a result it is an important ideal to look at and the algebraic invariants of it are of primary importance. There are only relatively few results on this ideal, apart from the results in the particular case of graph binomial ideals. The generalized binomial edge ideal is defined as the ideal generated by a collection of $2$-minors of a generic $2 \times n$ matrix. The determinantal ideal is of classical importance to algebraic geometry and commutative algebra. Apart from the geometric significance, another motivation to study these ideals comes from their connection to conditional independence ideals, a topic in algebraic statistics.
Let $m,n\geq 2$ be integers and $G$ be a simple graph on the vertex set $[n]$. Let $X=(x_{ij})$ be an $m\times n$ matrix of indeterminates and $S=K[X]$ be the polynomial ring in the variables $x_{ij}, 1\leq i \leq m, 1\leq j\leq n$, over a field $K$. For $1\leq k< l\leq m$ and $\{i,j\} \in E(G)$ with $1\leq i< j\leq n$ we set,
\begin{center}
    $p^{kl}_{ij}=[k,l][i,j]=x_{ki}x_{lj}-x_{li}x_{kj}$.
\end{center}
The ideal $\mathcal{J}_{G}=(p^{kl}_{ij}: 1\leq k <l \leq m, \{i,j\}\in E(G))$ is called the generalized binomial edge ideal.
In \cite{GBEI_Rauh}, it is proved that $\mathcal{J}_{G}$ is a radical ideal and its minimal primes are determined by the sets with the cut point property of $G$. Regarding the minimal free resolution of generalized binomial edge ideal, not much is known. But there have been a few key works in this area. Madani and Kiani in \cite{conjecture_pair}, proved that $\mathcal{J}_{G}$ has linear resolution if and only if $m=2$ and $G$ is a complete graph.
 In \cite{GBEI_Rida}, the authors proved the equality of the depth of generalized binomial edge ideal of generalized block graph and the depth of its initial ideal. Also they found the regularity of these two classes of generalized block graphs. Since a tree is a block graph, one can get an equivalent condition, which is given in \cite{GBEI_Rida}. Also, for a path, regularity of generalized binomial edge ideal is computed in \cite{GBEI_Rida}, and independently computed in \cite{GBEI_Arvind}. Regularity of the generalized binomial edge ideal for the complete graph is known in \cite{GBEI_Arvind}. In \cite{amata_generalized}, the authors have proved that the generalized binomial edge ideal is Cohen-Macaulay when the graph is a complete graph.
 
 Another class of binomial ideals related to the binomial edge ideal is the parity binomial edge ideal which was introduced in \cite{Kahle_parity}. The permanental binomial edge ideal and Lovasz-Saks-Schrjiver (in short, LSS) ideal in \cite{orthogonal_rep} are also important classes of binomial edge ideals associated to graphs. Not much work has been done on these classes of ideals and the authors of the present article feels that much work on these are also important to further understanding of the homological properties of the binomial ideals. The parity binomial edge ideal has generated some recent interest among the combinatorial commutative algebra community. Let $G$ be a simple graph with edge set $E(G)$ and $S=K[x_1,x_2,\ldots,x_n,y_1,y_2,\ldots,y_n]$ be a polynomial ring, the parity of binomial edge ideal of $S$ is denoted by $\mathcal{I}_G$ and is defined as $\mathcal{I}_G=(x_i x_j-y_i y_j: \{i,j\}\in E(G))$, whereas, the permanental edge ideal is denoted by $\Pi_G$ and is defined as $\Pi_G=(x_i y_j+x_jy_i: \{i,j\}\in E(G))$. Let $d\geq 1$ be an integer, then the Lovasz-Saks-Schrjiver ideal is denoted by $L_{G}^{K}(d)$ and is defined as,
\begin{center}
$L_{G}^{K}(d)=(\displaystyle \sum_{l=1}^{d} x_{il} x_{jl}: \{i,j\}\in E(G))$
\end{center}
in the polynomial ring $K[x_{kl}: 1\leq k\leq n, 1\leq l\leq d]$. For $d=1$, $L_G^{K}(d)$ it coincides with the edge ideal of graph $G$. For $d=2$, LSS ideal of $G$ is the binomial ideal defined as $L_G=(x_i x_j + y_i y_j: \{i,j\}\in E(G))\subset K[x_1,x_2,\ldots,x_n,y_1,y_2,\ldots,y_n]$. It is known that $L_G$ is a radical ideal if char$(K)\neq 2$, in \cite{orthogonal_rep}. Also, the authors in loc. cit. determined the primary decomposition of $L_G$ when $\sqrt{-1} \notin K$ and char$(K)\neq 2$. In \cite{Kahle_parity}, the authors proved that the parity binomial edge ideal is radical if and only if the graph is bipartite. Also, in the same paper they have determined the minimal primes and the primary decomposition. In \cite{LSS_Conca}, the authors proved that $L_G^K(2)$ is complete intersection if and only if $G$ does not contain a claw or an even cycle (Theorem 1.4). Also they have shown that $L_G^K(3)$ is prime if and only if $G$ does not contain a claw or $C_4$. In \cite{orthogonal_rep}, authors computed a Gr\"{o}bner basis of the permanental edge ideal of graphs and also they have shown that it is a radical ideal. In \cite{LSS_arvind}, the author characterizes the graphs whose parity binomial edge ideal is a complete intersection, see Theorem 3.2 and 3.5 of \cite{LSS_arvind}. In the loc. cit., the author characterized the graphs for LSS ideals, parity binomial edge ideal and permanental edge ideals that are complete intersections. In \cite{Arvind_parity_reg}, the authors have given a lower bound for the regularity of parity binomial edge ideal of graphs, and also classified the graphs whose parity binomial edge ideal has regularity 3. In the same article they have characterized the graphs whose parity binomial edge ideal have pure resolutions. In \cite{Hilbert_Parity}, the authors have given an explicit formula for the Hilbert-Poincare series of the parity binomial edge ideal of a complete graph.

\section{Preliminaries}
In this section, we recall basic definitions, notations, and terminologies of graphs and algebra which will be needed further.
\subsection{Basic notions from Algebra}
Let $S=K[x_1,x_2,\ldots,x_n]$ be a polynomial ring in $n$ variables over a field $K$. Then $S$ is called a graded ring, graded by degree, if $S=\oplus_{m\geq 0} S_m$, where $S_m$ is generated by all monomials of degree $m$, and if $I$ is a graded ideal of $S$, then $R=S/I$ is called graded algebra.

A graded free resolution of a finitely generated $S$-module $M$ is an exact sequence
\begin{center}
 $\mathbb{F}:  \ldots F_i  \overset{\phi_i}\longrightarrow F_{i-1} \longrightarrow \ldots \longrightarrow
F_1 \overset{\phi_1}\longrightarrow F_0 \overset{\phi_0}\longrightarrow M \longrightarrow 0$
\end{center}
where $F_i=\oplus_{j\in \mathbb{Z}} S(-j)^{\beta_{ij}}, \forall j$ are the free $S$-modules. The numbers $\beta_{ij}=\beta_{ij}(M)$ are called the graded Betti numbers of $M$, and $\beta_{i}=\sum_{j} \beta_{ij}$ is called the total $i$th Betti number of $M$.

The projective dimension of $M$ is defined by
\begin{center}
    $pd(M)=max \{i : \beta_{ij} \neq 0, \text{for some }   i\}$
\end{center}
The regularity of $M$ is defined by
\begin{center}
    $reg (M)=max \{j-i : \beta_{ij} \neq 0 \}$.
\end{center}
The depth of $M$ is defined by
\begin{center}
    $depth(M)= n- \text{proj dim} (M)$
\end{center}
If $\mathbf{m}=(x_1,x_2,\ldots,x_n)$ is the maximal ideal of $S$, then free resolution is called minimal if $\phi_{i+1}(F_{i+1})\subset \mathbf{m} F_i, \forall i$. Minimal free resolution of a finitely generated module is finite which follows from the following theorem.
\begin{theorem}[Hilbert Syzygy Theorem]
Let $S=K[x_1,x_2,\ldots,x_n]$ be a polynomial ring of $n$ variables over a field $K$. Then the graded minimal free resolution of a graded finitely generated $S$-module is finite, and its length is at most $n$.
\end{theorem}
A finitely generated $S$-module $M$ has a $d$-linear resolution if the minimal free resolution is of the form
\begin{center}
  $0 \rightarrow S(-d-i)^{\beta_i}\rightarrow \ldots \rightarrow S(-d-1)^{\beta_1}\rightarrow S(-d)^{\beta_0}\rightarrow M
 \rightarrow 0$
\end{center}
where the degree of each generators of $M$ is $d$.\\
The module $M$ has a pure resolution if the minimal free resolution is of the following form
\begin{center}
 $0 \rightarrow S(-d_i)^{\beta_i}\rightarrow \ldots \rightarrow S(-d_1)^{\beta_1}\rightarrow S(-d_0)^{\beta_0}\rightarrow M
 \rightarrow 0$
\end{center}
such that $0<d_0<d_1<\ldots<d_i$ are all integers.\\
An ideal $I$ is said to be Cohen-Macaulay if the Krull dimension coincide with the depth.

\subsection{Basic notions of Graph}
Let $G$ be a finite simple graph on $n$-vertices without isolated vertices. A subgraph $H$ of $G$ is called {\it induced}, if for every vertices $u$ and $v$ in $H$ whenever $\{u,v\}$ is an edge of $G$, it is an edge of $H$. A graph $G$ is called {\it chordal graph} if every induced subgraph of $G$ has cycle of length 3, and it is called {\it co-chordal} if its complement is chordal. A complete graph is called {\it clique}. The largest number of the clique is called {\it clique number} of $G$.
 A vertex $v$ of $G$ whose deletion from the graph gives a graph with more connected components, then $v$ is called a {\it cut point} of $G$. A subset $T\subset [n]$ is said to have a {\it cut point property} (in short, cut point set) for $G$, if for every 
$v \in T$, 
$c(T \setminus \{v\}) < c(T)$, where $c(T)$ is the number of connected components of the restriction of $G$ to $[n] \setminus T$. A cut set of a graph $G$ is a subset of the vertices whose deletion increase the number of connected components of $G$.

Let $G_1$ and $G_2$ be two graphs with vertex sets $V_1$ and $V_2$ and edge sets $E_1$ and $E_2$ respectively. The {\it join} of two graphs $G_1$ and $G_2$ denoted by $G_1*G_2$, is a graph defined on the vertex set $V_1\cup V_2$ and on the edge set 
    $E_1 \cup E_2\cup \{\{x,y\}: x\in V_1, y\in V_2\}$.
For a graph $G$, a maximal subgraph of $G$ without cut vertex is called a {\it block} of $G$. A graph $G$ is said to be a {\it block graph} if each block of $G$ is a clique. 
\begin{definition}(See, \cite{Arvind_Jahangir})
a graph $G$ is said to be a quasi-block graph if $G$ satisfies the following:
\begin{enumerate}
    \item Each block of $G$ is either a clique or a quasi-block
    \item If $v$ is an internal vertex of a quasi-block $B$, then for any $u\in N_G(v) \setminus V(B)$, $u$ is not an internal vertex for any block.
\end{enumerate}
\end{definition}
A graph $G$ is called an {\it interval graph} if every vertex $v\in V(G)$ can be labelled with a real closed interval $I_v=[a_v,b_v]$ in such a way that two distinct vertices $v, w \in V(G)$ are adjacent if their corresponding intervals have non-empty intersection. The vertices of the interval graph $G$, are intervals namely $V(G)=\{ I_1, \ldots, I_r\}$, where $I_j=[a_j,b_j]$ with $a_j\leq b_j$, for all $1\leq j \leq r$.

Let $\mathbb{N}_0=\mathbb{N}\cup\{0\}$. Also, let $J_0=[0], J_i=[i-1, i], i=1,\ldots,k$ and $I_j=[a_j,b_j]$,  such that $a_j \in \mathbb{N}_0$ and $I_j\subseteq [0,k]\ \{k\},$ for all $j=1,2,\ldots,r$. Then we call the interval graph on the vertex set $\{J_i\}_{i=0}^{k}\cup \{I_j\}_{j=1}^{r}$ a connected strong interval graph or simply a {\it strong interval graph}.

\section{Regularity of binomial edge ideal of graphs}
In this section, we will provide an overview of regularity of binomial edge ideal of graphs.
In \cite{binomial_reg2} and \cite{binomial_reg3}, authors have characterized regularity of binomial edge ideal of 2 and 3 respectively and its initial ideal. In the next theorem the authors characterize when binomial edge ideal and its initial ideal have a linear resolution that is reg$(J_G)=2$.
\begin{theorem}[\cite{sparse_ideal}, Theorem 1.4 and \cite{binomial_reg2}, Theorem 2.1]
Let $G$ be a simple graph, and $<$ be any term order. Then the following are equivalent.
\begin{enumerate}
   \item $J_G$ has a linear resolution.
    \item in$_<(J_G)$ has a linear resolution.
    \item $G$ is a complete graph.
\end{enumerate}
\end{theorem}
In \cite{CM_binomial}, authors conjectured that the regularity of binomial edge ideal of a graph is equal to the regularity of its initial ideal. This conjecture holds for the join of two graphs.
\begin{theorem}[\cite{binomial_reg3}, Theorem 2.1]\label{reg_join}
Let $G_1$ and $G_2$ be graphs on disjoint vertex sets $V_1$ and $V_2$ respectively, not both complete, let $<$ be any term order on ring $S$. Then
\begin{enumerate}
    \item reg$(J_{G_1*G_2})$=max$\{reg (J_{G_1}), reg (J_{G_2}), 3\}$
    \item reg $(in_< J_{G_1*G_2})$=max$\{reg(in_< J_{G_1}), reg (in_< J_{G_2}),3\}$
\end{enumerate}
\end{theorem}
As an application of Theorem \ref{reg_join}, they have characterized binomial edge ideal of regularity 3.
\begin{theorem}[\cite{binomial_reg3}, Theorem 3.2]
Let $G$ be a non-complete graph with $n$ vertices, and no isolated vertices. Then reg$(J_G)=3$ if and only if either
\begin{enumerate}
    \item $G=K_r \sqcup K_s$, with $r,s\geq 2$ and $r+s=n$, or
    \item $G=G_1*G_2$, where $G_i$ is a graph with $n_i<n$ vertices such that $n_1+n_2=n$ and reg$(J_{G_i})\leq 3$, for $i=1,2$.
\end{enumerate}
\end{theorem}
In \cite{conjecture_pair}, Madani-Kiani conjectured that the regularity of binomial edge ideal of a graph is at most $c(G)+1$, where $c(G)$ is the number of maximal cliques of $G$. This conjecture is proved for closed graphs, generalized block graphs, some class of chordal graphs and for fan graphs of complete graph.
Let $H$ be a connected closed graph on $[n]$ such that $S/J_{H}$ is Cohen-Macaulay. Then by (\cite{CM_binomial} Theorem 3.1), there exist integers $1=a_1<a_2\ldots<a_s<a_{s+1}=n$ such that $F_i=[a_i,a_{i+1}]$, for $1\leq i\leq s$ and $F_1,\ldots,F_s$ is a leaf order of $\triangle(H)$. Set $e=\{1,n\}$. The graph $G=H\cup \{e\}$ is called the quasi-cycle associated with $H$.

Let $H$ be a connected closed graph on the vertex set $[m]$ such that $S_{H}/J_{H}$ is Cohen-Macaulay. Then by (\cite{CM_binomial}, Theorem 3.1), there exist integer $1=a_1<a_2<\ldots<a_s<a_{s+1}=m$ such that $F_i=[a_i,a_{i+1}], 1\leq i\leq s$ and $F_1,F_2,\ldots, F_S$ is a leaf order of $\triangle(H)$. Let $F_{s+1}=[m,n]\cup\{1\}$. The graph on the vertex set $[n]$ and edge set $E(G)=E(H)\cup \{\{i,j\}: i\neq j, i,j \in F_{s+1}\}$ is called a semi-cycle graph associated with $H$.
\begin{definition}
A block $B$ of a graph $G$ is said to be a semi-block if $B$ is a semi-cycle with $B\neq K_3$. A graph $G$ is said to be a semi-block graph if all except one block are cliques and the block which is not a clique is a semi-block.
\end{definition}
In \cite{Arvind_Jahangir}, the author proved this conjecture for quasi block graphs and for semi block graphs.
\begin{theorem}[\cite{Arvind_Jahangir}, Theorem 3.7 and Theorem 3.11]
For a quasi-block graph as well as semi block graph $G$,  reg$(S/J_G)\leq c(G)$.
\end{theorem}
\begin{definition}
The Jahangir graph denoted by $J_{m,n}$ is a graph on the vertex set $[m,n+1], m\geq 1, n\geq 3$ such that the induced subgraph on $[mn]$ is $C_{mn}$ and the neighbourhood of the vertex $mn+1$ is $\{1,m+1,\ldots, m(n-1)+1\}$.
\end{definition}
In \cite{Arvind_Jahangir}, the author proved the same conjecture independently for chordal graph see, Theorem 3.15 and for Jahangir graphs see, Theorem 4.5. In \cite{construction_HibiBEI}, Hibi-Matsuda conjectured that for a graph $G$ on $n$-vertices, reg$(S/J_G)\leq deg h_{S/J_G}(t)$.  Characterization of generalized block graphs which admit unique extremal Betti number is given in \cite{reg_Gblock}. In \cite{Arvind_Jahangir}, author provide a sufficient condition for this conjecture.
\begin{theorem}[\cite{Arvind_Jahangir}, Theorem 5.1]
Let $G$ be a connected graph on vertex set $[n]$. If $S/J_G$ admits a unique extremal Betti number, then reg$(S/J_G)\leq deg h_{S/J_G}(t)$.
\end{theorem}
In \cite{Arvind_Jahangir}, Example 5.3 author has given a counter example for this conjecture. Also the author posed a question which is the following:
\begin{question}
When does the binomial edge ideal of a graph admit a unique extremal Betti number?
\end{question}
In \cite{reg_block}, authors proved that regularity of $S/J_G$ is less than or equal to $n-1$, for some classes of trees. An improved upper bound for the regularity of binomial edge ideal of trees is given as a corollary of the Theorem 0.1, \cite{reg_trees}
\begin{corollary}[\cite{reg_trees}, Corollary 0.1]
Let $T$ be a tree on $[n]$ with spine $P$ of length $l$. Let $e_2$ denote the number of edges that are not in $P$ and with both end points having degree at most 2 and $d_3$ denote the number of vertices, not in $P$, and having degree at least 3. Then
\begin{center}
    $reg(S/J_{T})\leq e_2 + l+ 2d_3$.
\end{center}
\end{corollary}
\begin{definition}
Let $I$ and $J$ be two proper ideals of a local regular ring $R$, they are called directly linked, we denote by $I\sim J$, if there exists a regular sequence $z=z_1,\ldots,z_n$ in $I\cap J$ such that $J=(z): I$ and $I=(z):J$. We say that $I$ and $J$ belong to the same linkage class if there exist a sequence a direct links 
\begin{center}
    $I=I_0 \sim I_1 \sim \ldots \sim I_m=J$.
\end{center}
$I$ is said to be linkage class complete intersection, in short licci, if $J$ is a complete intersection ideal.

\end{definition}
If $I$ and $J$ are linked, then if one of them is Cohen-Macaulay then other is too. In particular, any licci ideal is always Cohen-Macaulay. 
A necessary condition for a homogeneous ideal in a polynomial ring which is Cohen-Macaulay as well as licci is given in \cite{necc_licci}.
\begin{theorem}[\cite{necc_licci}, Corollary 5.13]
Let $I$ be a Cohen-Macaulay homogeneous ideal in a standard graded polynomial ring $S=K[x_1,\ldots, x_n]$ with the graded maximal ideal $\mathfrak{m}$. If $I_{\mathfrak{m}}\subset R=S_{\mathfrak{m}}$ is licci, then
\begin{center}
    $reg (S/I)\geq (height I-1)(indeg I-1)$
\end{center}
where $indeg I$ is the initial degree of the ideal $I$, that is $indeg I=min \{i: I_{i}\neq 0 \}$.
\end{theorem}
In \cite{reg_licci}, the authors have given a new upper bound for the regularity of $S/J_{G}$.
\begin{theorem}[\cite{reg_licci}, Theorem 2.1]
Let $G$ be a connected graph on vertex set $[n]$. Then 
\begin{center}
    $reg (S/J_{G})\leq n- dim \triangle (G)$
\end{center}
where $\triangle(G)$ is a simplicial complex of $G$.
\end{theorem}
For a disconnected graph $G$, the authors have found an upper bound, and it is the following
\begin{theorem}[\cite{reg_licci}, Corollary 2.2]
Let $G$ be a graph on $n$ vertices with the connected components $G_1,\ldots, G_c$. Then
\begin{center}
    $reg (S/J_{G})\leq n- (dim \triangle (G_1)+\ldots+dim \triangle(G_c))$.
\end{center}
\end{theorem}
Next theorem has improved the upper bound given by Matsuda and Murai in \cite{reg_length}.

\begin{figure}[H]
\begin{center}
\begin{tikzpicture}[every node/.style={circle, scale=.4}][font=\Large] 
\node [draw,shape=circle,scale=0.8][fill](one) at (1,0)  {};
\node [draw,shape=circle,scale=0.8][fill](two) at (2,0)  {};
\node [draw,shape=circle,scale=0.8][fill](three) at (3,0)  {};
\node [draw,shape=circle,scale=0.8][fill](four) at (4,0)  {};
\node [draw,shape=circle,scale=0.8][fill](five) at (0,2)  {};
\node [draw,shape=circle,scale=0.8][fill](six) at (0,3)  {};
\node [draw,shape=circle,scale=0.8][fill](seven) at (0,4)  {};
\node [draw,shape=circle,scale=0.8][fill](eight) at (0,5)  {};
\node [draw,shape=circle,scale=0.8][fill](nine) at (-1,0) {};
\node [draw,shape=circle,scale=0.8][fill](ten) at (-2,0)  {};
\node [draw,shape=circle,scale=0.8][fill](ele) at (-3,0)  {};
\node [draw,shape=circle,scale=0.8][fill](twe) at (-4,0)  {};
\draw (one) -- (two) node[midway,above] {$e_{1}$};
\draw (three) -- (four) node[midway,above] {$e_{r}$};
\draw (three)--(four);
\draw (five) -- (six) node[midway,right] {$e_{1}^{''}$};
\draw(seven) -- (eight) node[midway,right] {$e_{t}^{''}$};
\draw (nine) -- (ten) node[midway, above] {$e_{1}^{'}$};
\draw (ele) -- (twe) node[midway, above] {$e_{s}^{'}$};
\draw (one)--(five)--(nine)--(one);
\draw [dotted](two)--(three);
\draw [dotted](six)--(seven);
\draw [dotted](ten)--(ele);
\end{tikzpicture}
\caption{Licci graphs}
\label{One}
  \end{center}
  \end{figure}
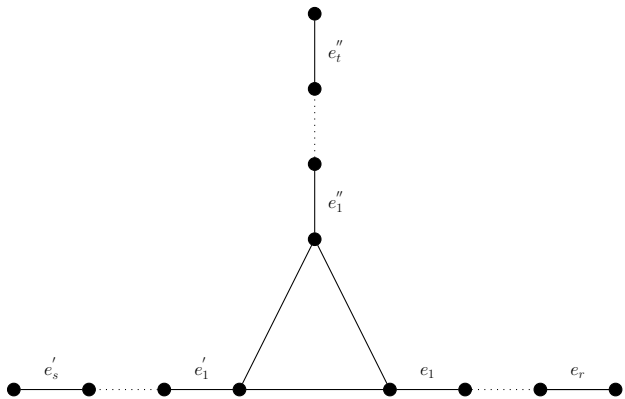

\begin{theorem}[\cite{reg_licci}, Theorem 3.5]\label{licci}
Let $G$ be a connected graph on the vertex set $[n]$. Then the following are equivalent
\begin{enumerate}
    \item $(J_G)_{\mathfrak{m}} \subset R=S_{\mathfrak{m}}$ is licci.
    \item $J_G$ is Cohen-Macaulay and $n-2 \leq reg (S/J_{G}) \leq n-1$
    \item $G$ is a path graph or it is isomorphic to one of the graphs in Figure \ref{One}, where $r, s, t$ are non-negative integers. In other words, $G$ is a triangle with possibly some paths connected to some of its vertices.
\end{enumerate}
\end{theorem}
An equivalent statement of \ref{licci} for chordal graph is the following.
\begin{theorem}[\cite{reg_licci}, Theorem 4.2]
Let $G$ be a connected chordal graph on the vertex set $[n]$. Then the following are equivalent
\begin{enumerate}
    \item $(J_{G})_{\mathfrak{m}}\subset R$ is licci.
    \item $J_G$ is Cohen-Macaulay and $n-2 \leq reg (S/J_{G}) \leq n-1$.
    \item $J_G$ is unmixed and $n-2 \leq reg (S/J_{G}) \leq n-1$.
    \item $G$ is a path graph or it is isomorphic to a graph in Figure \ref{One}.
\end{enumerate}
\end{theorem}
An explicit formula for the regularity of binomial edge ideal for the block graphs in terms of combinatorics of the graph is an open problem. An algorithm for dimension of binomial edge ideal for block graph is given in \cite{dimension}, Algorithm 2.5. 

For a flower graph, extremal Betti number is given in \cite{dimension}, Theorem 3.4. As a consequence, one can get regularity of binomial edge ideal for flower graph $F(v)$.
\begin{corollary}[\cite{dimension}, Corollary 3.5]
Let $F(v)$ be a flower graph, then
\begin{center}
    $reg S/J_{F(v)}=i(F(v))+cdeg(v)-1$
\end{center}
\end{corollary}
where $cdeg(v)$ is the clique degree of $v$, is the number of maximal cliques to which $v$ belongs and $i(F(v))$ is the inner vertices of $F(v)$.
\begin{definition}
Let $G$ be a block graph. If $G$ has no flower graph as an induced subgraph, then $G$ is called flower free graph.
\end{definition}
\begin{definition}
Let $G$ be a block graph and $F(v)$ be a flower graph as an induced subgraph of $G$. $F(v)$ is called an end-flower of $G$ if $G=G_1 \cup \ldots \cup G_c$, where $c=cdeg(v)$, such that $G_i \cap G_j=\{v\}$, for all $1\leq i <J\leq c$ and $G_2,\ldots, G_c$ are flower free graphs.
\end{definition}
The formula for the regularity of binomial edge ideal of block graphs.
\begin{theorem}[\cite{dimension}, Theorem 4.2]
Let $G$ be a block graph, $v_1, \ldots, v_r\in V(G)$
\begin{center}
    $H_j=G\ \{v_1,\ldots, v_j\}$
\end{center}
$j=1,2,\ldots,r$ and $H_0=G$. If
\begin{enumerate}
    \item $F(v)$ is an end-flower for $H_{j-1}$, for all $j=1,\ldots,r$
    \item $H_r$ is flower-free
\end{enumerate}
then
\begin{center}
$reg S/J_{G}=reg S/J_{H_{r}}=c+i(H_{r})$    
\end{center}
where $c$ is the number of connected components of $H_r$ which are not isolated vertices.
\end{theorem}
In \cite{reg_Gblock}, the author has given improved lower and upper bound for the regularity of binomial edge ideal of generalized block graphs. In the same article, they have also characterized generalized block graph in which binomial edge ideal admits a unique extremal Betti number.
\begin{theorem}[\cite{reg_Gblock}, Theorem 3.11]
Let $G$ be a connected indecomposable generalized block graphs. Then the following are equivalent
\begin{enumerate}
    \item $S/J_{G}$ admits a unique extremal Betti number.
    \item For any $v\in V(G)$, $F_{h,k}(v)$ is not an induced subgraph of $G$, for any $h,k \geq 0$ with $h+k \geq 3$.
\end{enumerate}
In this case, reg$(S/J_{G})=m(G)+1$.
\end{theorem}
where $F_{h,k}(v)$ is a flower graph obtained by gluing of $h$ copies of the complete graph $K_3$ and $k$ copies of the star graph $K_{1,3}$ at a common free vertex $v$ and $m(G)$ is the number of minimal cut sets of the graph $G$.

Let $deg_G(v)$ denote the degree of a vertex $v$ of a graph $G$. If $deg_G(v)=1$, then we say $v$ is a pendant vertex. Let for $v\in V(G)$, the number of maximal cliques of $G$ which contain $v$, be denoted by $cdeg_G(v)$ and the number of pendant vertices adjacent to $v$ be denoted by $pdeg_G(v)$. A vertex 
$v\in V(G)$ such that $cdeg_G(v)=pdeg_G(v)+1$ with 
$pdeg_G(v)\geq 1$ is said to be of type 1, and of type 2 if $cdeg_G(v)\geq pdeg_G(v)+2$.
 Let $\alpha(G)$ be the number of vertices of type 1, and $pv(G)$ be the number of pendant vertices of $G$.

With the same notations as in above, in \cite{reg_Gblock}, the author obtained upper bound for the regularity of binomial edge ideal of connected indecomposable generalized block graphs.
\begin{theorem}[\cite{reg_Gblock}, Theorem 4.5]
Let $G$ be a connected indecomposable generalized block graph on $[n]$ vertices, which is not a star graph. Then
\begin{center}
    reg$(S/J_{G})\leq c(G)+\alpha(G)-pv(G)$
\end{center}
where $c(G)$ is number of maximal cliques of $G$.
\end{theorem}
Let $G$ be a graph on the vertex set $V(G)$, and for subset $T\subseteq V(G)$ the induced subgraph of $G$ on the vertex set $V(G)\setminus T$, be denoted by $G-T$. Let $N_{G}(v)$ denote the set of neighbours of the vertex $v$ in $G$. We set $\hat{G}=G \setminus \mathcal{I}s(G)$, where $\mathcal{I}s(G)$ denotes the set of isolated vertices of $G$ and $K_n$ denote the complete graph on $n$ vertices. Associated to a vertex $v$ in $V(G)$, there is a graph denote by $G_v$, with the vertex set $V(G)$ and the edge set 
\begin{center}
    $E(G)\cup \{\{u,w\}: \{u,w\} \subseteq N_{G}(v)\}$
\end{center}
In \cite{Proof_conjecture}, the authors give a general upper bound for regularity of binomial edge ideal. For that they have introduced a map, called compatible map, and the definition is as the following:
\begin{definition}
Let $\mathcal{G}$ be the set of all graphs. We call a map $\phi: \mathcal{G} \rightarrow \mathbb{N}_0$, compatible, if it satisfies the following conditions
\begin{enumerate}
    \item $\phi(\hat{G})\leq \phi(G)$, for all $G \in \mathcal{G}$
    \item If $G=\sqcup_{i=1}^{t} K_{n_i}, n_i \geq 2$, for every $1\leq i\leq t$, then $\phi(G)\geq t$.
    \item If $G\neq \sqcup_{i=1}^{t} K_{n_i}$, then there exists $v\in V(G)$ such that
    \begin{enumerate}
        \item $\phi(G-v)\leq \phi(G)$ and
        \item $\phi(G_v)< \phi(G)$.
    \end{enumerate}
\end{enumerate}
\end{definition}
\begin{theorem}[\cite{Proof_conjecture}, Theorem 2.4]
Let $G$ be a graph on $[n]$ and $\phi$ a compatible map. Then
\begin{center}
    $reg (S/J_G)\leq \phi(G)$.
\end{center}
\end{theorem}
Next we will see a new upper bound for regularity of binomial edge ideal which is clique disjoint edge set in graph.
\begin{definition}\label{defn}
Let $G$ be a graph and $\mathcal{H}\subseteq E(G)$ with the property that no two element of $\mathcal{H}$ belong to a clique of $G$. Then we call the set $\mathcal{H}$, a clique disjoint edge set in $G$. we set
\begin{center}
    $\eta(G)=max\{ |\mathcal{H}|: \mathcal{H} \text{is a clique disjoint edge set in} $G$\}$
\end{center}
\end{definition}
\begin{theorem}[\cite{Proof_conjecture}, Corollary 2.7]
Let $G$ be a graph on vertex set $[n]$. Then 
\begin{center}
    $reg(S/J_{G})\leq \eta(G)$.
\end{center}
\end{theorem}
So, by above result we get a sharp bound for regularity of binomial edge ideal, as $\eta(G)\leq c(G)$, for every graph $G$.\\
In \cite{conjecture_chordal}, authors proved that regularity of strong interval graph of $G$ coincide with $\mathcal{L}(G)$ as well as $c(G)$, where $\mathcal{L}(G)$ is the sum of the lengths of the longest induced paths of connected components of $G$.
\begin{theorem}[\cite{conjecture_chordal}, Corollary 4.3]
Let $G$ be a chordal graph. Then the following are equivalent
\begin{enumerate}
    \item $reg(S/J_G)=\mathcal{L}(G)=c(G)$
    \item $G$ is a strong interval graph.
\end{enumerate}
\end{theorem}

\section{Depth of Binomial edge ideal}
In general, it is hard to compute another algebraic invariant such as the depth of binomial edge ideals. Regarding the depth of binomial edge ideal not much have been done. In the next result we will see the formula for the depth of cone graph. The definition of cone graph is the following
\begin{definition}
Let $H$ be a graph on vertex set $[n]$. The cone of $v$ on $H$ is denote by $v*H$, is a graph with vertex set $V(v*H)=V(H)\sqcup \{ v \}$ and edge set $E(v*H)=E(H)\sqcup \{\{ u,v\}: u \in V(H)\}$. 
\end{definition}
Let $G=v*H$ and $S_H=K[x_i,y_i:i\in V(H)]$ and $S=S_H[x_v,y_v]$. In \cite{depth_arvind}, authors has proved that under the process of taking cone on connected graph depth remains invariant .
\begin{theorem}[\cite{depth_arvind}, Theorem 3.4]
Let $H$ be a connected graph on the vertex set $[n]$ and let $G=v*H$ be the cone graph. Then 
\begin{center}
    $depth_S(S/J_{G})=depth_{S_{H}}(S_{H}/J_{H})$
\end{center}
In particular, if $H$ is Cohen-Macaulay, then $G$ is almost Cohen-Macaulay.
\end{theorem}
For the binomial edge ideal of cone graph, depth formula is the following.
\begin{theorem}[\cite{depth_arvind}, Theorem 3.9]
Let $G=v*H$, where $H$ is a disconnected graph on vertex set $[n]$. Then
\begin{center}
    $depth_S{S/J_{G}}=min\{depth_{S_{H}}(S_{H}/J_{H}), n+2\}$.
\end{center}
\end{theorem}
Depth formula of binomial edge ideal for join of two graphs can be found in \cite{depth_arvind}, Theorem 4.1, 4.3, 4.4. \\

 Let $H_{\mathfrak{m}} ^{i}(S/J_G)$ denote the $i$th local cohomology module of $S/J_G$, where maximal ideal associated to it is $\mathfrak{m}=(x_1,\ldots,x_n,y_1,\ldots,y_n)$. Then 
\begin{center}
 depth$(S/J_G)=min\{i: H_{\mathfrak{m}}^{i}(S/J_G)\neq 0\}$.  
 \end{center}
In \cite{Betti_cycle}, authors have shown that depth $(S/J_G)=n$, for $n>3$. In \cite{CM_binomial}, authors proved that depth $(S/J_{G})=n+1$, for a connected block graph $G$. Later in \cite{depth_GBG}, authors computed the depth of generalized block graphs. In \cite{depth_arvind}, authors have given the formula for depth of join product of two graphs $G_1$ and $G_2$, which is denoted by $G_1*G_2$. In \cite{graph_connectivity}, the authors has given an upper bound in terms of graphs invariant for the depth of binomial edge ideal of a graph. They proved that for a non-complete connected graph $G$, depth $(S/J_G)\leq n-\kappa(G)+2$, where $\kappa(G)$ denotes the vertex connectivity of $G$.

In \cite{small_depth}, authors proved Hochster type formula for the local cohomology modules of binomial edge ideal, which is based on \cite{Hochster_BEI},Theorem 3.9.
\begin{theorem}(\cite{small_depth}, Theorem 3.6)
Let $G$ be a graph on $[n]$, and $\mathcal{Q}(G)$ be the poset associated to $J_G$. Then we have the $\mathbb{K}$-isomorphism
\begin{center}
    $H_{\mathfrak{m}}^{i}(S/J_G)\cong \oplus _{q\in \mathcal{G}} H_{\mathfrak{m}}^{d_q} (S/q)^{\oplus M_{i,q}}$
\end{center}
where $d_q=dim S/q$ and $M_{i,q}=dim_{\mathbb{K}} H^{\tilde i-d_q-1}((q,1_{\mathcal{Q}_{G}}); \mathbb{K})$.
\end{theorem} 
In the next theorem, the  authors have given a lower bound for the depth of binomial edge ideal.
\begin{theorem}[\cite{small_depth}, Theorem 5.2]
Let $G$ be a graph on $[n]$. Then 
\begin{center}
    depth$(S/J_G)\geq 4r+\sum_{i=1}^{n} r_i(i+1)$,
\end{center}
where $r$ is the number of non-complete connected components of $G$, and $r_i$ is the number of complete connected components of $G$ of size $i$, for every $1\leq i\leq n$.
\end{theorem}
Next, the authors characterized the graphs for depth$(S/J_G)=4$.
\begin{theorem}[\cite{small_depth}, Theorem 5.3]
Let $G$ be a graph on $[n]$ with $n\geq 4$. Then the following are equivalent.
\begin{enumerate}
    \item depth$(S/J_G)=4$
    \item $G=G^{'}* 2K_1$, for some graph $G^{'}$,
\end{enumerate}
 where $2K_1$ is the graph consisting of two isolated vertices.
\end{theorem}

In \cite{graph_connectivity}, authors proved that for a connected non-completed graph $G$,
\begin{center}
    depth$(S/J_G)\leq n+2-\kappa(G)$,
\end{center}
where $\kappa(G)$ is the graph connectivity of $G$. In \cite{malayeri_depth}, authors proved that 
\begin{center}
    $f(G)+d(G)\leq depth (S/J_G)$.
\end{center}
In \cite{jayanthan_depth}, authors have characterized graphs for which depth of $S/J_G$ attain the lower bound.
\begin{theorem}[\cite{jayanthan_depth},Corollary 5.6]
Let $G$ be a unicyclic graph or a quasi-cycle graph such that depth$(S/J_G)=d(G)+f(G)$. Then either $G$ is chordal or $G$ has an induced $C_4$.
\end{theorem}
Next question arises is, whether there exists any graph $G$ such that $f(G)+d(G)=n+2-\kappa(G)$ or not. Existence of such graphs are characterized by Hibi and Saeedi Madani in \cite{hibi_diameter}. Also, in \cite{jayanthan_depth}, author characterizes connected non-complete graphs $G$ with the property that $f(G)+d(G)+1=n+2-\kappa(G)$.
\begin{theorem}[\cite{jayanthan_depth}, Theorem 3.9]
Let $G$ be a graph on $[n]$ such that $d(G)+f(G)+1=n+2-\kappa(G)$.
\begin{enumerate}
    \item If $\kappa(G)=1$, then either $G$ is chordal or $G$ has precisely one induced $C_4$ and has no induced $C_l$, for $l\geq 5$.
    \item If $\kappa(G)\geq 2$ and $d(G)=2$, then $G$ is a chordal graph.
    \item If $\kappa(G)=2$ and $d(G)=3$, then either $G$ is chordal or $G$ has precisely one induced $C_4$, and has no induced $C_l$, for $l\geq 5$.
\end{enumerate}
\end{theorem}
The authors in (\cite{jayanthan_depth}, Section 4) described the structure of the graphs with the property $f(G)+d(G)+1=n+2-\kappa(G)$ and calculated the depth of $S/J_G$. Also in the same paper authors posed the following question.
\begin{question}[\cite{jayanthan_depth}, Question 5.5]
If $G$ is a graph containing an induced cycle of length at least 5, then is depth$(S/J_G)\geq d(G)+f(G+1)$?
\end{question}

\section{Results on the other algebraic invariants}
\subsection{Betti number}
Betti number is one of the homological invariant which will help us to know about its structure. Researchers have computed Betti number of binomial edge ideals. Zofar and Zahid calculated for cycles in \cite{Betti_cycle}, Schenzel and Zafar calculated for complete bipartite graphs in \cite{Betti_bipartite} and Jayanthan et. al calculated for trees and unicycle graphs in \cite{Betti_tree}. In \cite{depth_arvind} authors have calculated Betti number for cone graph [see, Theorem 3.10]. As a consequence of this, they have calculated Betti number for wheel graph [see, \cite{depth_arvind}, Corollary 3.11].

In \cite{Betti_tree}, authors computed the first Betti number of binomial edge ideal of trees (see, Theorem 3.1) and second Betti number of unicyclic graphs (see, Theorem 3.4). In the same paper, they have calculated first syzygy of binomial edge ideal for cycle $C_n$, Theorem 3.5. Next, they have calculated first syzygy of binomial edge ideal for unicyclic graphs (Theorem 3.6, 3.7) and for trees.

Let $A\subset [n]$, and $i\in A$, define $P_A(i)=|\{j\in A: j\leq i\}|$. Here, $P_A$ denotes the position of an element in $A$ when the elements are arranged in the ascending order. Now we are ready to state the following theorem.
\begin{theorem}[\cite{Betti_tree}]
Let $G$ be a tree on $n$-vertices. Then the first syzygy of $J_G$ is minimally generated by elements of the form
\begin{enumerate}
    \item $f_{i,j}e_{\{k,l\}}-f_{k,l}e_{\{i,j\}}$, where $\{i,j\}, \{k,l\}\in E(G)$ and $\{e_{\{i,j\}}: \{i,j\}\in E(G)\}$ is the standard basis of $S(-2)^{n-1}$.
\item $(-1)^{P_A(j)} f_{k,l} e_{\{i,j\}}+ (-1)^{P_A(k)} f_{j,l} e_{\{i,k\}}+(-1)^{P_A(l)} f_{j,k} e_{\{i,l\}}$, where $A=\{i,j,k,l\}\in \mathcal{C}_G$ with center at $i$.
\end{enumerate}
\end{theorem}
where $\mathcal{C}_G=\sum_{v\in V(G)} \binom{deg_G(v)}{3}$.

 Rees algebra of an ideal $I$ is $\mathcal{R}(I)=\oplus_{n\geq 0} I^n t^n$. The generators of the defining ideal of the Rees algebra of a graph is described by Villarreal in \cite{Rees}.
  He proved that $I(G)$ is of linear type, that is the Rees algebra is isomorphic to the symmetric algebra if and only if $G$ is either a tree or an odd unicyclic graph. But almost nothing is known about the Rees algebra of binomial edge ideal of a graph. It is known that for a connected graph $G, J_G$ is complete intersection if and only if $G$ is a path \cite{CM_binomial}.

Let $G$ be a graph on $n$-vertices and $J_G$ be its binomial edge ideal. Let $R=S[T_{\{i,j\}}\in E(G), i<j]$. Let $\delta : R \rightarrow S[t]$ be the $S$-algebra homomorphism given by $\delta(T_{\{i,j\}})=f_{i,j}t$. Then $Im(\delta)=\mathcal{R}(J_G)$ and $ker(\delta)$ is called the defining ideal of $\mathcal{R}(J_G)$.

The following theorem characterizes the trees whose binomial edge ideal are almost complete intersection.
\begin{theorem}[\cite{small_depth}, Theorem 4.3]
If $G$ is a tree which is not a path, then $J_G$ is an almost complete intersection ideal if and only if $G$ is obtained by adding an edge between two vertices of two paths.
\end{theorem}
In general, they have characterized for connected graph which is not a tree in Theorem 4.4. Also they have shown that the associated graded ring and the Rees algebra of almost complete intersections binomial edge ideal are Cohen-Macaulay.
\begin{theorem}[\cite{small_depth}, Theorem 4.7]
If $G$ is a graph such that $J_G$ is an almost complete intersection ideal, then $gr_S(J_G)$ and $\mathcal{R}(J_G)$ are Cohen-Macaulay.
\end{theorem}
Since, complete intersections are of linear type, binomial edge ideal of paths are linear type. In \cite{small_depth}, Proposition 4.9, authors proved that binomial edge ideal of $K_{1,n}$ is of linear type. In the same paper authors proved that in the polynomial ring over an infinite field, almost complete intersection homogeneous ideals are generated by $d$-sequence, Proposition 4.10. As a corollary of the above result we have the following.
\begin{corollary}[\cite{small_depth}, Corollary 4.11]
Let $G$ be a graph on $n$-vertices. If $J_G$ is an almost complete intersection ideal, then $J_G$ is generated by a $d$-sequence. In particular, $J_G$ is of linear type.
\end{corollary}
Also, they have obtained the defining ideal of Rees algebra of binomial edge ideal of cycles as a corollary, see [\cite{small_depth}, Corollary 4.13]. Authors have asked the following question:
\begin{question}[\cite{small_depth}, Question 4.16]
Classify all bipartite graphs whose binomial edge ideal are of linear type.
\end{question}
The above question is not true for tree. So, one can prove or disprove the following,
\begin{conj}[\cite{small_depth}, Conjecture 4.17]
If $G$ is a tree or a unicyclic graph, then $J_G$ is of linear type.
\end{conj}

\subsection{Construction of a graph}
In \cite{construction_Hibi}, Hibi et. al. constructed a graph $G$ such that for $1\leq b \leq r$, the regularity of monomial edge ideal of $G$ is $r$, and the number of its extremal Betti number is $b$. Also in \cite{construction_HibiBEI}, for a given pair $(r,s)$ with $1\leq r\leq s$, Hibi and Matsuda constructed a graph $G$ such that reg$(S/J_{G})=r$ and the degree of the $h$-polynomial of $S/J_{G}$ is $s$. In \cite{depth_arvind}, authors constructed the following.
\begin{theorem}[\cite{depth_arvind}, Theorem 5.4]
Let $r$ and $b$ be two positive integers with $1\leq b\leq r-1$. Then there exist a graph $G=G_{r,b}$ such that reg$(S/J_{G})=r$, and the number of extremal Betti number of $S/J_{G}$ is $b$.
\end{theorem}
\begin{question}
Does there exists a graph $G$ such that the projective dimension is bounded by a linear function of $b$ and $r$, where $r=reg(S/J_{G})$, and $b$ is the number of extremal Betti number of $S/J_{G}$?
\end{question}

\subsection{Cohen-Macaulay Binomial edge ideal}
Cohen-Macaulay binomial edge ideals are studied by many authors  \cite{CM_binomial}, \cite{construction_CM}. Also, Cohen Macaulay binomial edge ideal for bipartite graphs and block graphs were studied in \cite{CM_binomial}, \cite{graph_connectivity}, \cite{CM_Bipartite}, \cite{depth_GBG}.
But the classification of it in terms of underlying graph is not known. In \cite{CM_deviation}, authors have classified Cohen Macaulay and unmixed binomial edge ideal with deviation. As an extension of these results, in \cite{cactus}, authors classified Cohen-Macaulay and unmixed binomial edge ideal $J(G)$, where $G$ is a cactus graph i.e. a graph whose blocks are cycles.

In \cite{CM_binomial}, Ene, Herzog and Hibi proved that if $G$ is a closed graph then $S/J_{G}$ is Gorenstein if and only if $G$ is a path.
Motivated by this result, author in \cite{Gorenstein} obtained the similar result for connected graph $G$.
\begin{theorem}[\cite{Gorenstein}, Theorem 5.3]
Let $S=K[x_1,\ldots,x_n,y_1,\ldots,y_n]$. Suppose the char$(K)=p>0$. Let $G$ be a connected graph such that $S/J_{G}$ is Gorenstein. Then $G$ is a path.
\end{theorem}
The above result also holds for char$(K)=0$ (\cite{Gorenstein}, Theorem 5.4).
\subsection{Hilbert series}
The Hilbert series of binomial edge ideal of cycles and quasi cycles is computed in \cite{HS_quasi}, \cite{Betti_cycle}. Authors in \cite{HS_Rinaldo} have computed the Hilbert-Poincare series of the binomial edge ideal of some Cohen-Macaulay bipartite graphs, Fan graphs (see Proposition 3.8). In \cite{HS_join}, authors obtained the Hilbert series for decomposable graphs in terms of indecomposable components.
\begin{definition}
A graph $G$ is said to be decomposable, if there exist induced subgraphs $G_1$ and $G_2$, such that $G=G_1 \cup G_2, V(G_1)\cap V(G_2)=\{v\}$, and $v$ is a free vertex of $G_1$ and $G_2$.
\end{definition}
\begin{theorem}[\cite{HS_join}, Theorem 3.2]\label{HS_deomposable}
Let $G=G_1 \cup G_2$ be a graph on vertex set $[n]$. Then 
\begin{center}
    $Hilb_{S/J_G}(t)=(1-t)^2 Hilb_{S_1/ J_{G_1}}(t) Hilb_{S_2/J_{G_2}(t)}$,
    \end{center}
    where $S_i=K[x_j,y_j: j\in V(G_i)]$, for $i=1,2$
\end{theorem}
They also have obtained dimension and multiplicity of $S/J_G$. As a corollary of \ref{HS_deomposable}, they have computed Hilbert series for a connected graph. In the same paper, authors computed Hilbert series for join of two graphs.
\begin{theorem}[\cite{HS_join}, Theorem 4.2]
Let $H$ and $H^{'}$ be two disconnected graphs on vertex sets $[p]$ and $[q]$, respectively. Let $G=H*H^{'}$ be the join of $H$ and $H^{'}$. Then
\begin{center}
  $Hilb_{S/J_{G}}(t)=Hilb_{S_H/J_H} (t)+Hilb_{S_{H^{'}}/J_{H^{'}}}(t)+ \frac{(p+r-1)t+1}{(1-t)^{p+q+1}}-\frac {(p-1)t+1}{(1-t)^{p+1}}-\frac {(q-1)t+1}{(1-t)^{q+1}}$
\end{center}
\end{theorem}
\subsection{Symbolic power}
In general, symbolic power and ordinary power of an ideal do not coincide. But for some classes of homogeneous ideals in polynomial rings it coincide, e.g. for a bipartite graph symbolic power and ordinary powers of an edge ideal coincide. In \cite{symbolic}, Ene and Herzog proved that for any binomial edge ideal with quadratic Gr\"{o}bner basis, the symbolic powers and ordinary powers of $J_G$ coincide.
Authors have shown that based on the condition of in$_<(J_{G})$ these two powers coincide.
\begin{theorem}[\cite{symbolic}, Theorem 3.3]
Let $G$ be a connected graph on the vertex set $[n]$. If $in_<(J_{G})$ is a normally torsion-free ideal, then 
\begin{center}
    $J_{G}^{(k)}=J_{G}^k$, for $k\geq 1$
\end{center}
\end{theorem}
As a consequence of the above result one can get the equality of these two powers for a closed graph, see, \cite{symbolic}, Corollary 3.4.
\subsection{Hochster Type formula}
Hochster formula originally appeared in \cite{Hochster_formula}, it provides a decompostion of the local cohomology modules $H_{\mathfrak{m}}^r (A/I)$ of the Stanley-Reisner ring $R/I$, associated to a squarefree monomial ideal $I\subseteq A=K[x_1,\ldots,x_n]$. It also describes the Hilbert series of the local cohomology modules of $A/I$. In \cite{Hochster_BEI} authors have given a Hochster type decomposition for the local cohomology modules associated to binomial edge ideal and is as the following.
\begin{theorem}[\cite{Hochster_BEI}, Theorem 3.9]
Let $A=K[x_1,\ldots,x_n, y_1,\ldots, y_n]$ be a polynomial ring over a field $K$. Let $J_G$ be the binomial edge ideal associated to a graph $G$ on the vertex set $[n]$. Let $Q_{J_{G}}$ be the poset associated to a minimal primary decomposition of $J_G$. Then the local cohomology modules with respect to $\mathfrak{m}$ of $A/J_{G}$ admit the following decompostion as $K$-vector spaces
\begin{center}
    $H_{\mathfrak{m}}^r (A/J_{G})\cong \oplus_{q \in Q_{J_{G}} H_{\mathfrak{m}}^{d_{q}}}(A/I_q)^{\oplus M_{r,q}}$
\end{center}
where $M_{r,q}=dim _{k} H^{\tilde r-d_{q}-1} ((q,1_{Q_{J_{G}}}); K)$. Moreover, we have a decomposition as graded $K$-vector spaces.
\end{theorem}
In \cite{Cartwright}, when authors comparing invariant of Cartwright-Strumfels ideals with its generic initial ideal they got the following
\begin{conj}[\cite{Cartwright}, Conjecture 1.14]\label{con}
Let $I\subseteq A$ be a $\mathbb{Z}^m$-graded Cartwright-Strumfels ideal and gin $(I)$ its $\mathbb{Z}^m$-graded generic initial ideal, with $m\leq n$. Then one has
\begin{center}
    $dim_K H_{\mathfrak{m}}^r (A/I)_a=dim_K H_{\mathfrak{m}}^r (A/gin (I))_a$.
\end{center}
for every $r\in \mathbb{N}$ and every $a\in \mathbb{Z}^m$.
\end{conj}
In \cite{Hochster_BEI}, the author has proved the conjecture \ref{con} for binomial edge ideal.
\begin{theorem}[\cite{Hochster_BEI}, Theorem 4.5]
Let $A=K[x_1,\ldots,x_n,y_1,\ldots,y_n]$ be a polynomial ring over $K$ and $\mathfrak{m}$ be its homogeneous maximal ideal. Let $J_G$ be binomial edge ideal associated to a graph $G$ on the vertex set $[n]$. Then 
\begin{center}
   $dim_K H_{\mathfrak{m}}^r (A/J_{G})_a=dim_K H_{\mathfrak{m}}^r (A/gin (J_{G}))_a$,
\end{center}
for every $r\in \mathbb{N}$ and every $a\in \mathbb{Z}^n$.
\end{theorem}
\section{Generalized Binomial edge ideal}
Generalized binomial edge ideal is introduced by Rauh in \cite{GBEI_Rauh}. It is the ideal generated by a collection of $2$-minors in generic matrix. Motivation to study these ideals comes from their connection to conditional independence ideals.
\begin{definition}
Let $m,n\geq 2$ be integers. Let $G$ be a simple graph on the vertex set $[n]$. Let $X=(x_{ij})$ be an $m\times n$ matrix of indeterminate and denote $S=K[X]$ the polynomial ring in the variables $x_{ij}, 1\leq i \leq m, 1\leq j\leq n$, over a field $K$. For $1\leq k< l\leq m$ and $\{i,j\} \in E(G)$ with $1\leq i< J\leq n$ we set
\begin{center}
    $p^{kl}_{ij}=[k,l][i,j]=x_{ki}x_{lj}-x_{li}x_{kj}$
\end{center}
The ideal $\mathcal{J}_{G}=(p^{kl}_{ij}: 1\leq k <l \leq m, \{i,j\}\in E(G)$ is called the generalized binomial edge ideal.
\end{definition}
In \cite{GBEI_Rauh}, it was proved that $\mathcal{J}_{G}$ is a radical ideal and its minimal primes are determined by the sets with the cut point property of $G$. Not much is known about minimal free resolution of generalized binomial edge ideal. Madani and Kiani in \cite{conjecture_pair}, proved that $\mathcal{J}_{G}$ has linear resolution if and only if $m=2$ and $G$ is a complete graph.

Let $G$ be a generalized block graph. Let $A_i(G)$ be the collection of cut sets of $G$ of cardinality $i$, $i=1,2,\ldots,\omega(G)-1.$ We denote $a_i(G)=|A_i(G)|$. In \cite{GBEI_Rida}, authors proved the equality of depth of generalized binomial edge ideal of generalized block graph and depth of its initial ideal, and also found the regularity of these two generalized block graph.
\begin{theorem}[\cite{GBEI_Rida}, Theorem 3.3]\label{reg}
Let $m,n\geq 2$, and let $G$ be a connected generalized block graph on the vertex set $[n]$. Then
\begin{enumerate}
    \item depth $(S/{\mathcal{J}_G})$=depth $(S/{in_<{(\mathcal{J}_G)}})=n+(m-1)-\displaystyle \sum_{i=2}^{\omega(G)-1} a_i (G)$
    \item If $m\geq n$, then reg $ S/{\mathcal{J}_G}$=reg $ S/{in_< (\mathcal{J}_G)}$
    \item If $m<n$, then reg $S/{\mathcal{J}_G}\leq$ reg $S/{in_<{\mathcal{J}_G}}\leq n-1$.
\end{enumerate}
\end{theorem}
As a consequence, one can get similar statements for block graph, as we have for block graph $a_i(G)=0,$ for all $i>1$. Since tree is a block graph, we can have equivalent result of the theorem \ref{reg} for tree. For a path graph regularity of generalized binomial edge ideal is the following 
\begin{theorem}[\cite{GBEI_Rida}, Corollary 3.6]
Let $m,n\geq 2$. If $G$ is a path graph on the vertex set $[n]$, then 
\begin{center}
    reg $S/{\mathcal{J}_G}$=reg $S/{in_<{(\mathcal{J}_G)}}=n-1$.
\end{center}
\end{theorem}
Independently, in \cite{GBEI_Arvind} author has given regularity bound for generalized binomial edge ideal. The author has defined in the following way

Let $G_1$ and $G_2$ be graphs on the vertex set $[m]$ and $[n]$ respectively. Let $e=\{i,j\}\in E(G_1)$, and $e^{'}=\{k,l\}\in E(G_2)$. Assign $2$-minor $p_{e,e^{'}}=x_{i,k}x_{j,l}-x_{i,l}x_{j,l}$ to the pair $(e,e^{'})$. The binomial edge ideal of the pair $(G_1, G_2)$ is, 
\begin{center}
$J_{G_1,G_2}=(p_{e,e^{'}}: e\in E(G_1), e^{'} \in E(G_2))$.    \end{center}
The generalized binomial edge ideal of a graph $G$ is the binomial edge ideal of the pair $(K_m,G)$. When $m=2$, the ideal $J_{K_2, G}=J_G$ is the classical binomial edge ideal of $G$.

In \cite{conjecture_pair}, Madani and Kiani conjectured that, if $G$ is a graph on vertex set $[n]$, then reg$(S/{J_{K_{m},G}})\leq $ min $\{ \binom {m}{2} c(G), e(G)\}$, where $c(G)$ denote the number of maximal cliques of $G$, and $e(G)$ denote the number of edges in $G$. They have proved the conjecture for closed graphs.
In \cite{GBEI_Arvind}, author proved this conjecture for chordal graphs.
\begin{theorem}[\cite{GBEI_Arvind}, Theorem 3.13]
Let $G$ be a connected chordal graph on the vertex set $[n]$. If $m<n$, then
\begin{center}
    $reg(S/J_{K_{m},G})\leq min \{ (m-1) c(G), n-1\}$.
\end{center}
\end{theorem}
In \cite{GBEI_Arvind}, author computed regularity of generalized binomial edge ideal for complete graph $K_n$.
\begin{proposition}[\cite{GBEI_Arvind}, Proposition 3.3]
Let $G=K_n$. Then reg $(S/J_{K_{m},K_n})=min \{ m-1, n-1\}$.
\end{proposition}
Also, for a connected graph $G$, in the same paper author obtained an upper bound for the regularity of generalized binomial edge ideal.
\begin{theorem}[\cite{GBEI_Arvind}, Theorem 3.6]
Let $G$ be a connected graph on vertex set $[n]$ and $m\geq 2$. Then 
\begin{center}
    $reg (S/{J_{K_{m},G}})\leq n-1$.
\end{center}
\end{theorem}
For chordal graph, the upper bound for regularity is the following.
\begin{theorem}[\cite{GBEI_Arvind}, Theorem 3.11]
Let $G$ be a connected chordal graph on the vertex set $[n]$. Then reg$(S/J_{K_{m}, G})\leq (m-1) c(G)$.
\end{theorem}

In \cite{GBEI_Arvind}, author conjectured for generalized binomial edge ideal for a connected graph, for $m<n$.
\begin{conj}[\cite{GBEI_Arvind}, Conjecture 3.15]
Let $G$ be a connected graph of vertex set $[n]$. Then
\begin{center}
    $reg(S/J_{K_{m},G})\leq min \{ (m-1) c(G), n-1\}$, if $m<n$.
\end{center}
\end{conj}
One can ask the following question for pair of a graphs $G_1$ and $G_2$.
\begin{question}[\cite{GBEI_Arvind}, Question 3.16]
Let $G_1$ and $G_2$ be graphs on the vertex set $[m]$ and $[n]$ respectively. Then
\begin{center}
    $reg(S/ J_{G_1,G_2})\leq min \{ (m-1) c(G_2), (n-1) c(G_1)\}$.
\end{center}
\end{question}
 In \cite{GBEI_Rida},\cite{Hibi_pair},  for some special cases unmixedness and Cohen-Macaulayness of binomial edge ideal of a pair of graphs are characterized. It is well-known that all Cohen-Macaulay ideal are unmixed. In \cite{amata_generalized}, author characterizes all the unmixed ideal $J_{G,m}$. Before that we will see the definition of power-cycle. 

Let $S\subseteq \{1,2,\ldots,\lfloor \frac{n}{2} \rfloor\}$. The circulant graph $G=C_n(S)$ is a simple graph with $V(G)=\mathbb{Z}_n=\{1,2,\ldots,n-1\}$ and $E(G=\{\{i,j\}:|j-i|_n \in S\})$, where $|k|_n=min \{|k|, n-|k|\}$. $C_n(S)$ is the circulant graph of order $n$ with generating set $S$ and $|k|_n$ is the circular distance modulo $n$. It is known that if the generating set $S=\{\pm 1,\pm 2,\ldots,\pm d\}$, where $1\leq d\leq \lfloor \frac{n}{2} \rfloor$ is a given integer, then the circulant graph $C_n(S)$ is equivalent to the $d$-th power of $C_n$, where two vertices are adjacent if and only if their distance is at most $d$. In this case, we will denote $C_n(S)$ by $C_n(1,2,\ldots,d)$, and is called $d$th power-cycle.
\begin{theorem}[\cite{amata_generalized}, Theorem 3.14]
Let $G=C_n(1,2,\ldots,r)$ be a non-complete power-cycle. Then $J_{G,m}$ is unmixed if and only if $m$ is odd, $n\in \{m+1,\ldots,\frac{3m+1}{2}\}$ and $r=\frac{m-1}{2}$.
\end{theorem}
In the same paper the author characterizes when $J_{G,m}$ is Cohen-Macaulay.
\begin{theorem}[\cite{amata_generalized}, Corollary 4.3]
Let $J_{G,m}$ be a generalized binomial edge ideal with $m\geq 3$. Then $J_{G,m}$ is Cohen-Macaulay if and only if $G$ is the complete graph.
\end{theorem}
\section{Other Binomial Edge Ideals}
In this section, we will give an overview of other binomial ideals namely the parity binomial edge ideal, the permanental edge ideal and the Lovasz-Saks-Schrjiver ideal.
In \cite{Kahle_parity}, the authors have given Gr\"{o}bner basis for parity binomial edge ideal $\mathcal{I}_G$, see Theorem 3.6 and minimal primes in Theorem 4.15. Also they have proved when $\mathcal{I}_G$ is radical ideal.
\begin{theorem}[\cite{Kahle_parity}, Theorem 5.5] Let $G$ be a graph. If char$(K)\neq 2$, then $\mathcal{I}_G$ is a radical ideal.
\end{theorem}
In the same article they have calculated primary decomposition of parity binomial edge ideal. In \cite{LSS_arvind}, the author has given an alternative form of the results from \cite{LSS_Conca}, for LSS ideals in terms of parity binomial edge ideal. 
\begin{theorem}[\cite{LSS_arvind}, Theorem 3.2] Let $G$ be a bipartite graph on $[n]$. Then $L_G$ is complete intersection, if and only if $\mathcal{I}_G$ is complete intersection, if and only if $G$ is a disjoint union of paths.
\end{theorem}
They have also characterized the non-bipartite graphs whose LSS ideals, permanental edge ideal and parity binomial edge ideals are complete intersection.
\begin{theorem}[\cite{LSS_arvind}, Theorem 3.5] Let $G$ be a non-bipartite graph on $[n]$. Then $L_G$ is complete intersection, if and only if $G$ is an odd cycle, if and only if $\mathcal{I}_G$ is complete intersection.
\end{theorem}
In the same article, they have computed projective dimension of the parity binomial edge ideal of an odd unicyclic graphs.

\begin{theorem}[\cite{LSS_arvind}, Theorem 4.4] For a connected odd unicyclic graph on $[n]$, projective dimension of $\mathcal{I}_G$ is pd$(S/\mathcal{I}_G)=n$.
\end{theorem}
Next, they have computed second Betti number and first syzygy of LSS ideals of trees.
\begin{theorem}[\cite{LSS_arvind}, Theorem 5.1] Let $G$ be a tree on $[n]$. Then
\begin{center}
$\beta_2(S/L_G)=\beta_{2,4} (S/L_G)=\binom{n-1}{2}+ \displaystyle \sum_{v \in V(G)} \binom {deg_G(v)}{3}$.
\end{center}
\end{theorem}
\begin{theorem}[\cite{LSS_arvind}, Theorem 5.2]
Let $G$ be a tree on $[n]$. If $\{ e_{\{i,j\}}: \{i,j\}\in E(G)\}$ is standard basis for $S^{n-1}$, then the first syzygy of $L_G$ is generated by the following elements
\begin{enumerate}
    \item $g_{i,j} e_{\{k,l\}}-g_{k,l}e_{\{i,j\}}$, where $\{i,j\}\neq \{k,l\} \in E(G)$ and
    \item $(-1)^{p_A(j)} f_{k,l} e_{\{i,j\}}+(-1)^{p_A(k)} f_{j,l} e_{\{i,k\}}+(-1)^{p_A(l)} f_{j,k} e_{\{i,l\}}$,\\ where $A=\{i,j,k,l\}\in \mathcal{C}_G$, with center at $i$.
\end{enumerate}
\end{theorem}
Also, they have computed second Betti number and first syzygy of LSS ideals for odd unicyclic graphs in Theorem 5.3 and 5.4.
In \cite{Arvind_parity_reg}, the author has given upper bound for the regularity of parity binomial edge ideal of a graph and for an odd cycle.
\begin{theorem}[\cite{Arvind_parity_reg}, Theorem 3.2] Let $G$ be a connected bipartite graph on $n$-vertices. Then reg$(S/ \mathcal{I}_G)\leq n-1$. Moreover, reg$(S/ \mathcal{I}_G)=n-1$, if and only if $G$ is a path graph.
\end{theorem}
\begin{theorem}[\cite{Arvind_parity_reg}, Theorem 3.4] Let $G=C_n$, where $n$ is odd. Then reg$(S/ \mathcal{I}_G)=n$.
\end{theorem}
In the next theorem, the author has characterized the graphs whose parity binomial edge ideal have regularity three.
\begin{theorem}[\cite{Arvind_parity_reg}, Theorem 3.9] Let $G$ be a graph on $n$-vertices with no isolated vertex. Then reg$(S/ \mathcal{I}_G)=2$ if and only if either $G= K_2 \sqcup K_2$ or $G$ is a complete bipartite graph other than $K_2$.
\end{theorem}
In the following theorem, the author has characterizes graphs whose parity binomial edge ideal have pure resolution.
\begin{theorem}[\cite{Arvind_parity_reg}, Theorem 3.14] Let $G$ be a graph on $n$-vertices. Then $S/ \mathcal{I}_G$ has pure resolution if and only if $G$ is one of the following
\begin{enumerate}
    \item $G$ is a complete bipartite graph.
    \item $G$ is a disjoint union of some odd cycles and some paths.
\end{enumerate}
\end{theorem}
Let $H_{S/ \mathcal{I}_G}$ be the Hilbert function of $S/ \mathcal{I}_G$. Then the Hilbert-Poincare series of the $S$-module $S/ \mathcal{I}_G$ is
\begin{center}
    $HP_{S/ \mathcal{I}_G}(t)=\displaystyle \sum_{i \geq 0} H_{S/ \mathcal{I}_G}(i) t^i$.
\end{center}
From (\cite{Peeva_syzygy}, Theorem 16.2), the series has the following expression
\begin{center}
$HP_{S/ \mathcal{I}_G}(t)=\frac{P_{S/ \mathcal{I}_G}(t)}{(1-t)^n}$.
\end{center}
$P_{S/ \mathcal{I}_G}(t)$ is the Hilbert-Poincare polynomial of $S/ \mathcal{I}_G$ and has the following form
\begin{center}
    $P_{S/ \mathcal{I}_G}(t)= \displaystyle \sum_{i=0}^{p} \sum_{j=0}^{p+r} (-1)^i \beta_{i,j} (S/ \mathcal{I}_G) t^j$.
\end{center}
In the following theorem the authors have calculated Hilbert-Poincare polynomial of $S/ \mathcal{I}_{K_n}$.
\begin{theorem}[\cite{Hilbert_Parity}, Theorem 3.6] The Hilbert-Poincare polynomial of $S/ \mathcal{I}_{K_n}$ is
\begin{center}
$P_{S/ \mathcal{I}_{K_n}}(t)=2(1-t)^n + [-1+3t+ (\frac{n^2+n-6}{2})t^2+(\frac{n^2-3n+2}{2})t^3] (1-t)^{2n-3}$.
\end{center}
In particular, depth$(S/ \mathcal{I}_{K_n})\geq 3$ and reg$(S/ \mathcal{I}_{K_n})\leq 3$.
\end{theorem}

{\bf Acknowledgement}: The author is not funded by any authorities.

\bibliographystyle{abbrv}
\bibliography{main}
\end{document}